\documentclass[graybox]{svmult}

\usepackage{mathptmx}       
\usepackage{helvet}         
\usepackage{courier}        
\usepackage{type1cm}        
														
\usepackage{makeidx}         
\usepackage{graphicx}        
\usepackage{multicol}        
\usepackage[bottom]{footmisc}

\usepackage{newtxtext}       %
\usepackage{newtxmath}       

\usepackage[utf8]{inputenc}
\usepackage{xcolor}
\usepackage[english]{babel}
\usepackage{amsfonts,amsmath,amssymb,amstext}
\usepackage{mathtools}
\usepackage[tight,nice]{units}
\usepackage{url}
\usepackage{graphicx}
\usepackage[compatibility=false]{caption}
\usepackage{subcaption}
\usepackage{cancel}
\usepackage{hyperref}

\newtheorem{thm}{Theorem}
\numberwithin{thm}{section}

\newcommand{\pwm}{f}
\newcommand{\numteeth}{m}

\newcommand{\bu}{\mbox{\boldmath$u$}}
\newcommand{\bv}{\mbox{\boldmath$v$}}
\newcommand{\be}{\mbox{\boldmath$e$}}
\newcommand{\bff}{\mbox{\boldmath$f$}}

\newcommand{\bU}{\mbox{\boldmath$U$}}
\newcommand{\stabFnc}{\mbox{$\mathcal{R}$}}

\makeindex 

\begin{document}
	\title*{{A New Parareal Algorithm for Time-Periodic Problems with Discontinuous Inputs}}
	\titlerunning{Parareal for Time-Periodic Problems with Discontinuous Inputs} 
	\author{Martin J. Gander, Iryna Kulchytska-Ruchka, and Sebastian Sch{\"o}ps}
	\institute{Martin J. Gander \at Section de Math\'{e}matiques, University of Geneva, 2-4 Rue du Li\`{e}vre,  CH-1211 Geneva, Switzerland, \email{martin.gander@unige.ch}	
	\and Iryna Kulchytska-Ruchka and Sebastian Sch{\"o}ps \at Institut f\"ur Theorie Elektromagnetischer Felder, 
		Technische Universit{\"a}t Darmstadt, Schlossgartenstrasse 8, D-64289 Darmstadt, Germany, \email{kulchytska@gsc.tu-darmstadt.de} and \email{schoeps@temf.tu-darmstadt.de}
		}
	%
	%
\maketitle

\abstract*{The Parareal algorithm, which is related to multiple
    shooting, was introduced for solving evolution problems in a
  time-parallel manner. The algorithm was then extended to solve
    time-periodic problems. We are interested here in time-periodic
  systems which are forced with quickly-switching discontinuous
  inputs. Such situations occur, e.g., in power engineering when
  electric devices are excited with a pulse-width-modulated signal. In
  order to solve those problems numerically we consider a recently
  introduced modified Parareal method with reduced coarse
  dynamics. Its main idea is to use a low-frequency smooth input for
  the coarse problem, which can be obtained, e.g., from Fourier
  analysis. Based on this approach, we present and analyze a new
  Parareal algorithm for time-periodic problems with
  highly-oscillatory discontinuous sources. We illustrate the
    performance of the method via its application to the simulation
  of an induction machine.  }

\section{Introduction} \label{section:introduction}

{Time-periodic problems appear naturally in engineering
  applications}. For instance, the {time-periodic} steady-state
behavior of an electromagnetic device is often the main interest in
electrical engineering, because devices are operated most of their
life-time in this state. Depending on the size and complexity of the
underlying system, {the} search for {a time-}periodic solution
might {however} be prohibitively expensive. {Special techniques}
were developed for {the} efficient computation of {such
  solutions, like} the \textit{time-periodic explicit error correction
  method} \cite{Katagiri_2011aa}{, which} accelerates calculations
by correcting the solution after each half period{, or the method}
presented in \cite{Bermudez_2013aa}, {which} leads to faster
{computations} of periodic solutions by determining suitable
initial conditions.

The Parareal algorithm {was invented in \cite{Lions_2001aa} for the
  parallelization of evolution problems in the time direction.}  A
detailed {convergence analysis when applied to} linear ordinary and
partial differential equations {with smooth right-hand sides can be
  found} in~\cite{Gander_2007a}{, for nonlinear problems, see}
\cite{Gander_2008a}. {In \cite{Kulchytska-Ruchka_2018ac}, a new
  Parareal algorithm was introduced and analyzed} for problems with
discontinuous sources. The main idea of the method is to {use} a
smooth approximation of the original signal as the input {for the
  coarse propagator}. In \cite{Gander_2013ab}, a Parareal algorithm
for nonlinear time-periodic problems was presented and analyzed.
{Our interest here is in time-periodic steady-state solutions of
  problems with quickly-switching discontinuous excitation, for which
  we will introduce and study a new {periodic} Parareal algorithm.}
	
\section{Parareal for time-periodic problems with discontinuous inputs} 
\label{section:parareal_for_TP_pbms}

{We} consider a time-periodic problem {given by} a system of
ordinary differential equations (ODEs) of the form
\begin{equation} \label{eq:ode}
  \bu^{\prime}(t) = \bff(t, \bu(t)), \quad t \in \mathcal{I},\qquad
  \bu(0) = \bu(T),
\end{equation}
with the right-hand side (RHS) $\bff: \bar{\mathcal{I}} \times
\mathbb{R}^n \to \mathbb{R}^n$ such that $\bff(0, \bv)=\bff(T, \bv)$
for all $n$-dimensional vectors $\bv$ and the solution $\bu:
\bar{\mathcal{I}} \to \mathbb{R}^n$ on the time interval $\mathcal{I}
:= (0, T)$.

In power engineering, electrical devices are often excited with a
pulse-width-modulated (PWM) signal \cite{Bose_2006aa}{, which} is a
discontinuous function {with} quick{ly}-switching
dynamics. {F}or electromagnetic applications such as motors or
transformers a couple of tens of kHz might be used as the switching
frequency \cite{Niyomsatian_2017aa}. {To solve time-periodic
  problems of the form \eqref{eq:ode}} supplied with such inputs
{with our new {periodic} Parareal algorithm, we assume that} the
RHS {can be split into} a sufficiently smooth bounded function
$\bar{\bff}$ and the corresponding discontinuous remainder
$\tilde{\bff}$ as
\begin{equation}\label{eq:rhs_splitting}
  \bff(t, \bu(t)) := \bar{\bff}(t, \bu(t)) + \tilde{\bff}(t),\quad t\in\mathcal{I}.
\end{equation}
{We decompose} $[0,T]$ into $N$ subintervals $[T_{n-1},T_n],$
$n=1,\dots,N$ with $T_0=0$ and $T_N=T$, {and introduce} the fine
propagator $\mathcal{F}\big(T_n, T_{n-1},\mathbf{U}^{(k)}_{n-1}\big)$
{which computes an accurate solution {at time $T_n$} of the
  initial-value problem (IVP)}
\begin{align}\label{eq:ode_nf}
  \bu_n^{\prime}(t) = \bff(t, \bu_n(t)), \quad t \in (T_{n-1}, T_{n}], 
    \qquad\bu_n(T_{n-1}) = \bU^{(k)}_{n-1}.
\end{align}
The {corresponding} coarse propagator
$\bar{\mathcal{G}}\big({T_n}, T_{n-1},\mathbf{U}^{(k)}_{n-1}\big)$
computes an {inexpensive approximation of the solution at time
  $T_n$ of the} corresponding IVP having the reduced smooth RHS
$\bar{\bff}(t, \bu(t))$,
\begin{equation}\label{eq:ode_nc}
  \bar{\bu}_n^{\prime}(t) = \bar{\bff}(t, \bar{\bu}_n(t)),
    \quad t \in (T_{n-1}, T_{n}], 
    \qquad\bar{\bu}_n(T_{n-1}) = \bU^{(k)}_{n-1}.
\end{equation}
{Our new periodic Parareal algorithm then computes f}or
$k=0,1,\dots$ and $\ n=1,\dots,N$
\begin{align}
\mathbf{U}_0^{(k+1)}&=\mathbf{U}_N^{(k)},\label{eq:PPIC_init}\\
\mathbf{U}_n^{(k+1)}&=\mathcal{F}\big(T_n, T_{n-1},\mathbf{U}^{(k)}_{n-1}\big)+
\bar{\mathcal{G}}\big(T_n, T_{n-1},{\mathbf{U}^{(k+1)}_{n-1}}\big)-
\bar{\mathcal{G}}\big(T_n, T_{n-1},\mathbf{U}^{(k)}_{n-1}\big),\label{eq:PPIC_iter}
\end{align}
until the jumps at the synchronization points $T_n,$ $n=1,\dots,N-1$
as well as the periodicity error between $\mathbf{U}_0^{(k)}$ and
$\mathbf{U}_N^{(k)}$ are reduced {to} a given tolerance.  With an
initial guess $\mathbf{U}_0^{(0)}$ at $T_0${, the initial guess}
$\mathbf{U}_n^{(0)},$ $n=1,\dots,N$ {for the algorithm
  \eqref{eq:PPIC_init}-\eqref{eq:PPIC_iter}} could be calculated
using, for instance, the reduced coarse propagator $\bar{\mathcal{G}}$
as
\begin{equation}\label{eq:PPIC_approx_kinit}
  \bU_n^{(0)}:=\bar{\mathcal{G}}\big(T_n, T_{n-1},{\bU^{(0)}_{n-1}}\big),
    \quad n=1,\dots,N.
\end{equation}
We note that the correction \eqref{eq:PPIC_init} does not impose a
strict periodicity, but a relaxed one, since the end value at the
$k$th iteration $\mathbf{U}_N^{(k)}$ is used to update the initial
approximation $\mathbf{U}_0^{(k+1)}$ at the next iteration. This
approach was introduced for time-periodic problems in
\cite{Gander_2013ab} and was named PP-IC (which stands for
\textit{periodic Parareal with initial-value coarse problem}). In
contrast to this method, where both coarse and fine propagators solve
the IVP \eqref{eq:ode_nf}, our iteration
\eqref{eq:PPIC_init}-\eqref{eq:PPIC_iter} uses a reduced dynamics on
the coarse level, described by \eqref{eq:ode_nc}. Convergence of the
PP-IC algorithm was analysed in \cite{Gander_2013ab}. {We extend
  this analysis now to the new} Parareal iteration
\eqref{eq:PPIC_init}-\eqref{eq:PPIC_iter}, applied to a
one-dimensional model problem.

\section{Convergence of the new periodic Parareal iteration}\label{section:convergence_1D}

We consider the linear time-periodic scalar ODE
\begin{equation} \label{eq:ode_1D}
  u^{\prime}(t)+\kappa u(t) = f(t), \quad t \in (0,T), \qquad u(0) = u(T),
\end{equation}
with a $T$-periodic discontinuous RHS $f: [0,T] \to \mathbb{R}$, a
constant $\kappa\in\mathbb{R}:$ $\kappa>0$, and the {solution} function
$u: [0,T] \to \mathbb{R}$ {we want to compute}.

In order to investigate {the} convergence of the {new periodic}
Parareal {algorithm} \eqref{eq:PPIC_init}-\eqref{eq:PPIC_iter}
{applied to \eqref{eq:ode_1D}}, we introduce several
assumptions. Let the time interval $[0,T]$ be decomposed into
subintervals of {equal length} $\Delta T=T/N$. {We assume} that
the fine propagator $\mathcal{F}$ is exact, {and we can thus} write
the solution of the IVP for {the ODE in} \eqref{eq:ode_1D} at
$T_n$, starting from the initial {value} $U_{n-1}^{(k)}$ at
$T_{n-1}$ as
\begin{equation}\label{eq:soln_fine}
\mathcal{F}\big(T_n, T_{n-1},{U}^{(k)}_{n-1}\big)=e^{-\kappa\Delta T}{U}^{(k)}_{n-1}+\int\limits_{T_{n-1}}^{T_n}e^{-\kappa\left(T_n-s\right)}f(s)ds.
\end{equation}
{Next}, introducing a smooth {and} slowly-varying RHS $\bar{f}$
{by} $f=\bar{f}+\tilde{f}$, we let the coarse propagator
$\bar{\mathcal{G}}$ be a one-step method, applied to
\begin{equation}\label{eq:ode_nc_1D}
  \bar{u}_n^{\prime}(t)+\kappa\bar{u}(t) = \bar{f}(t),
  \quad t \in (T_{n-1}, T_{n}], 
  \qquad \bar{u}_n(T_{n-1}) = U^{(k)}_{n-1}.
\end{equation}
{Using the stability function $\stabFnc\left(\kappa\Delta T\right)$ of
  the one-step method, o}ne can then write
\begin{equation}
\label{eq:soln_coarse}
\bar{\mathcal{G}}\big(T_n, T_{n-1},{U}^{(k)}_{n-1}\big)=\stabFnc\left(\kappa\Delta T\right){U}^{(k)}_{n-1}+\xi\left(\bar{f},\kappa\Delta T\right),
\end{equation}
{where the} function $\xi$ corresponds to the RHS discretized with
the {one-step method}.  We also assume that
\begin{equation}
\label{eq:cond_coarse}
\left|\stabFnc\left(\kappa\Delta T\right)\right|+\left|e^{-\kappa\Delta T}-\stabFnc\left(\kappa\Delta T\right)\right|<1.
\end{equation}
{Using} \eqref{eq:soln_fine} and \eqref{eq:soln_coarse} {and
  following \cite{Gander_2013ab},} the errors
{$e^{(k+1)}_n:=u(T_n)-U^{(k+1)}_n$} of the {new periodic
  Parareal algorithm} \eqref{eq:PPIC_init}-\eqref{eq:PPIC_iter}
{applied to the model problem \eqref{eq:ode_1D} satisfy {for
    $n=1,2,\dots,N$} the relation}
\begin{align}\nonumber
e^{(k+1)}_n&=u(T_n)-\mathcal{F}\big(T_n, T_{n-1},{U}^{(k)}_{n-1}\big)
-\bar{\mathcal{G}}\big(T_n, T_{n-1},{U}^{(k+1)}_{n-1}\big)+\bar{\mathcal{G}}\big(T_n, T_{n-1},{U}^{(k)}_{n-1}\big)\\ \nonumber
&=e^{-\kappa\Delta T}u(T_{n-1})+\int_{T_{n-1}}^{T_n}e^{-\kappa\left(T_n-s\right)}f(s)ds-e^{-\kappa\Delta T}{U}^{(k)}_{n-1}-\int_{T_{n-1}}^{T_n}e^{-\kappa\left(T_n-s\right)}f(s)ds\\ \nonumber
&\quad-\left(\stabFnc\left(\kappa\Delta T\right){U}^{(k+1)}_{n-1}+\xi\left(\bar{f},\kappa\Delta T\right)\right)
+\left(\stabFnc\left(\kappa\Delta T\right){U}^{(k)}_{n-1}+\xi\left(\bar{f},\kappa\Delta T\right)\right)\\
&=\stabFnc\left(\kappa\Delta T\right)e^{(k+1)}_{n-1}+\left(e^{-\kappa\Delta T}-\stabFnc\left(\kappa\Delta T\right)\right)e^{(k)}_{n-1}.
\label{eq:e_n_k_plus_one}
\end{align}
  Similarly, the initial error satisfies $e^{(k+1)}_0=e^{(k)}_N$. {A
  key observation here is} that there is no explicit reference to the
  right-hand sides $f$ or $\bar{f}$ in \eqref{eq:e_n_k_plus_one}{:
  the corresponding terms cancel both between the exact solution and
  the (exact) fine solver, and also between the two coarse solvers!}  
	{Collecting the errors in the error vector
  $\be^{(k)}:=\left(e^{(k)}_0,e^{(k)}_1,\dots,e^{(k)}_N\right)^T$, we
  obtain from \eqref{eq:e_n_k_plus_one}} the {same} fixed-point
iteration {as in \cite{Gander_2013ab}, namely}
\begin{equation}
  \label{eq:errorIter}
	\be^{(k+1)}= S\be^{(k)},
\end{equation}
{where} the matrix $S$ is given by
\begin{equation}\label{eq:S}
S=
\begin{bmatrix}
1 &   &  & 0 \\
-\stabFnc\left(\kappa\Delta T\right)  & 1 &  &  \\
 &  \ddots   &  \ddots & \\
 &    &  -\stabFnc\left(\kappa\Delta T\right)    & 1
\end{bmatrix}^{-1}
\begin{bmatrix}
0 &   &   & 1 \\
e^{-\kappa\Delta T}-\stabFnc\left(\kappa\Delta T\right)  & 0 &  &  \\
 &  \ddots   &  \ddots & \\
 &    &  e^{-\kappa\Delta T}-\stabFnc\left(\kappa\Delta T\right)    & 0
\end{bmatrix}.
\end{equation}
{T}he asymptotic convergence factor {of the fixed-point
  iteration \eqref{eq:errorIter} describing our new periodic Parareal
  algorithm \eqref{eq:PPIC_init}-\eqref{eq:PPIC_iter} applied to the
  periodic problem \eqref{eq:ode_1D} is therefore given by}
\begin{equation}\label{eq:asym_conv_factor}
  \rho_{\mathrm{asym}}(S)=\lim\limits_{k\to\infty}\left({\|\be^{(k)}\|}/{\|\be^{(0)}\|}\right)^{1/k}.
\end{equation}
\begin{thm}[{Convergence estimate of the new periodic Parareal algorithm}]\label{thm:PPIC}
{Let $[0,T]$ be partitioned into $N$ equal time intervals with}
$\Delta T=T/N$. {Assume the} fine propagator {to be exact as} in
\eqref{eq:soln_fine}{, and t}he coarse {propagator} to be a one-step
method {as in} \eqref{eq:soln_coarse} {satisfying}
\eqref{eq:cond_coarse}.  Then the asymptotic convergence factor
\eqref{eq:asym_conv_factor} of the new {periodic} Parareal
algorithm \eqref{eq:PPIC_init}-\eqref{eq:PPIC_iter} is bounded {for all $l\geq1$} {by}
\begin{equation} \label{eq:conv_estimate}
\rho_{\mathrm{asym}}(S)<x_l,\,\ \mathrm{with}\,\ x_l=\left(\left|\stabFnc\left(\kappa\Delta T\right)\right|x_{l-1}+
\left|e^{-\kappa\Delta T}-\stabFnc\left(\kappa\Delta T\right)\right|\right)^{\frac{N}{N+1}}\ {\mathrm{and}\,\ x_0=1.}
\end{equation}
\begin{proof}
  Since the errors {of the new periodic Parareal algorithm} satisfy the
  {same} relation \eqref{eq:errorIter} {as in \cite{Gander_2013ab}}, 
	the proof follows {by the same arguments as in}  \cite{Gander_2013ab}.
\end{proof}
\end{thm}
We note that under the assumption \eqref{eq:cond_coarse}{, the}
operator $S$ is a contraction \cite{Gander_2013ab}, which ensures
convergence of the {new periodic Parareal algorithm} \eqref{eq:PPIC_init}-\eqref{eq:PPIC_iter}.

\section{Numerical experiments for a model problem}\label{section:example_1D}

\begin{figure}[t]
  \begin{subfigure}[t]{.475\textwidth}
		\centering
		\includegraphics[height=4.7cm]{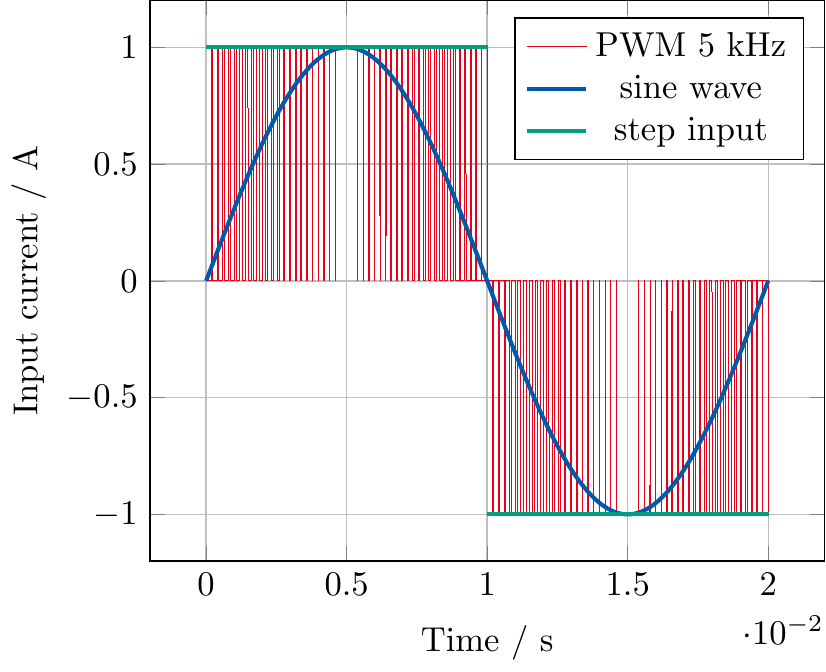}
	\end{subfigure}%
	\hspace*{0.04\textwidth}
	\begin{subfigure}[b]{.475\textwidth}
		\centering
		\includegraphics[height=4.7cm]{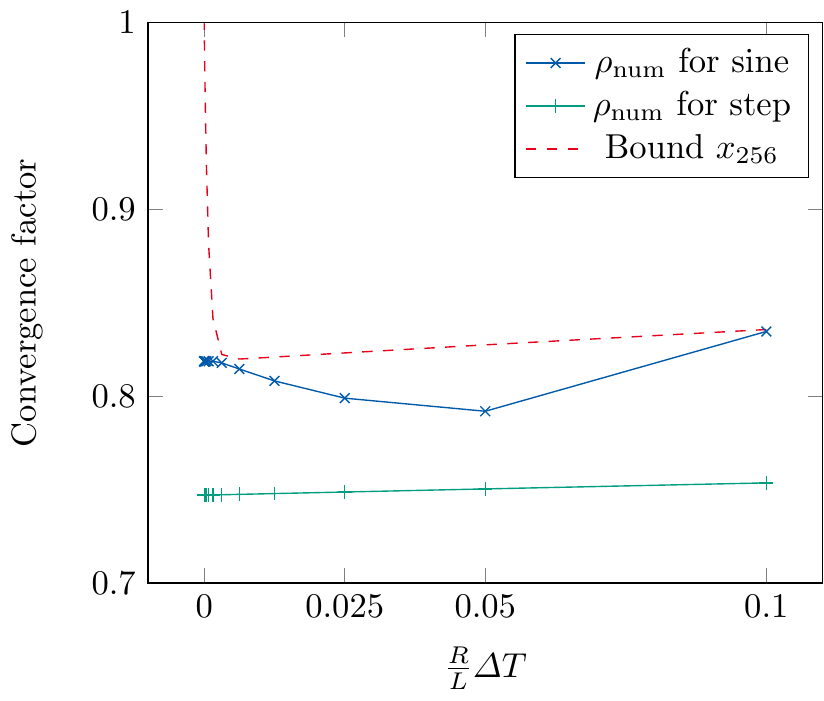}
	\end{subfigure}
  \caption{Left: PWM excitation \eqref{input_pwm} of $5$ kHz $(m=100)$
    and two coarse inputs \eqref{input_sine}, \eqref{input_step}.
    Right: convergence {factor} of the {new} periodic Parareal
    algorithm \eqref{eq:PPIC_init}-\eqref{eq:PPIC_iter} with
    reduced coarse dynamics: sinusoidal waveform \eqref{input_sine}
    and step function \eqref{input_step}{, together with our
      theoretical bound}.  }
	\label{fig:conv_pwm_sine_step}
\end{figure}

In this section we {illustrate our} convergence theory {for the new periodic Parareal algorithm with} a
periodic problem {given by} an RL-circuit model{, namely}
\begin{equation} \label{RL_circuit_model}
  R^{-1}\phi^\prime(t)+L^{-1}\phi(t)= \pwm_{\numteeth}\left(t\right),
  \quad t\in(0,T),
  \qquad \phi(0) = \phi(T),
\end{equation}
with the resistance $R=0.01$\;$\Omega$, inductance $L=0.001$\;H,
period $T=0.02$\;s, and $\pwm_{\numteeth}$ denoting the supplied PWM
current source (in A) of $20$ kHz, defined by
\begin{equation}\label{input_pwm}
\pwm_{\numteeth}(t)=\begin{cases}
\mathrm{sign}\left[\sin\left(\dfrac{2\pi}{T} t\right)\right],\ & s_m(t)-\left|\sin\left(\dfrac{2\pi}{T}t\right)\right|<0,\\
0,\ & \mathrm{otherwise,}
\end{cases}
\end{equation}
where $s_m(t)=\dfrac{\numteeth}{T}t -
\left\lfloor\dfrac{\numteeth}{T}t\right\rfloor,$ $t\in[0,T]$ is the
common sawtooth pattern with $m=400$ teeth. An example of the PWM
signal of $5$ kHz is shown in Figure~\ref{fig:conv_pwm_sine_step} on
the left.
This figure also illustrates the following two choices for
the coarse excitation:
\begin{align}\label{input_sine}
\bar{f}_\mathrm{sine}(t)&=\sin\left(\frac{2\pi}{T}t\right),\quad t\in[0,T],\\
\bar{f}_\mathrm{step}(t)&=
\begin{cases}
1,\ &t\in[0,T/2],\\
-1,\ &t\in(T/2,T].
\end{cases} \label{input_step}
\end{align}
We note that the step function \eqref{input_step} is discontinuous 
only at $t=T/2$. {This} does not lead to any difficulties, since {we}
use {Backward} Euler for {the time} discretization and we choose the
discontinuity to be located exactly at a synchronization point.

The coarse propagator $\bar{\mathcal{G}}$ then solves an IVP for the
equation $R^{-1}\phi^\prime(t)+L^{-1}\phi(t)=\bar{f}(t),\ t\in(0,T]$,
where the RHS $\bar{f}$ is one of the functions {in}
\eqref{input_sine} or \eqref{input_step}.  We {illustrate} the
estimate \eqref{eq:conv_estimate} by calculating the numerical
convergence factor
$\rho_{\mathrm{num}}:=\left(\|\be^{(K)}\|/\|\be^{(0)}\|\right)^{1/K}$
with the $l^\infty$-norm of the error $\be^{(k)}$ at iteration
$k\in\{0,K\}$ defined as \vspace{-0.075cm}
\begin{equation}
	\|\be^{(k)}\|=\max\limits_{0\leq n\leq
  N}|\phi(T_n)-\Phi_n^{(k)}|.
\end{equation}
\vspace{-0.075cm}\noindent
Here $\phi$ denotes the {time-periodic} 
steady-state solution of \eqref{RL_circuit_model} having the
same accuracy as the fine propagator, and $\Phi_n^{(K)}$ is the
solution obtained at the $K$th iteration when
\eqref{eq:PPIC_init}-\eqref{eq:PPIC_iter} converged up to a prescribed
tolerance. {The} stability function used for $x_l$ in 
\eqref{eq:conv_estimate} in case of Backward Euler is 
$\stabFnc\left(\frac{R}{L}\Delta T\right)=\left(1+\frac{R}{L}\Delta T\right)^{-1}.$

{On the right in} Figure~\ref{fig:conv_pwm_sine_step}{, we show
  the measured convergence factor} of the {new} Parareal iteration
\eqref{eq:PPIC_init}-\eqref{eq:PPIC_iter} for the two choices of the
coarse excitation \eqref{input_sine} and \eqref{input_step}. The fine
step size is chosen to be $\delta T=T/2^{18}\sim 7.63^{-8}$, while the
coarse step varies as $\Delta T=T/2^p$, $p=1,2,\dots,17$. We {also
show on the right in Figure~\ref{fig:conv_pwm_sine_step} the value
of} $x_{256}$ to be the bound in \eqref{eq:conv_estimate}. 
The graphs show that the theoretical estimate is indeed an upper bound 
for the numerical convergence factor for both coarse inputs (sine and step). 
However, one can observe that $x_{256}$ gives a sharper estimate in the 
case of the sinusoidal RHS \eqref{input_sine}, compared to the one 
defined in \eqref{input_step}. We also noticed that the number of 
iterations required till convergence of \eqref{eq:PPIC_init}-\eqref{eq:PPIC_iter} 
was the same (9 iterations on average for the values of $\Delta T$ considered) 
for both choices of the coarse input, while the initial error $\|\be^{(0)}\|$ 
was bigger with {the step coarse input} \eqref{input_step} than with the sinusoidal 
waveform \eqref{input_sine}. This led to a slightly smaller convergence factor
in case of the step coarse input due to the definition of $\rho_{\mathrm{num}}$.

\section{{Numerical experiments} for an induction machine}\label{section:example_machine}

{We now test the} performance of {our new periodic Parareal
  algorithm} with reduced coarse dynamics for {the}
simulation of a four-pole squirrel-cage induction motor, excited by a
three-phase PWM voltage source switching at $20$ kHz.  The model of
this induction machine {was} introduced in~\cite{Gyselinck_2001aa}.
We consider the no-load condition, when the motor operates with
synchronous speed.

{The} spatial discretization of the two-dimensional cross-section
of the machine with $n=4400$ degrees of freedom leads to {a}
time-periodic problem {represented by the} system of
differential-algebraic equations (DAEs)
\begin{align}\label{eq:dae}
	\mathbf{M}\mathrm{d}_t\mathbf{u}(t)
	+
	\mathbf{K}\left(\mathbf{u}(t)\right)\mathbf{u}(t)
	&=
	\mathbf{f}(t),\quad {t\in(0,T)},\\
	\mathbf{u}(0)&=\mathbf{u}(T),\label{eq:dae_TP}
\end{align}
with unknown $\mathbf{u}:[0,T]\to\mathbb{R}^{n},$ (singular) mass
matrix $\mathbf{M}\in\mathbb{R}^{n\times n},$ nonlinear stiffness
matrix
$\mathbf{K}\left(\cdot\right):\mathbb{R}^{n}\to\mathbb{R}^{n\times
  n},$ and the $T$-periodic RHS $\mathbf{f}:[0,T]\to\mathbb{R}^{n},$
$T=0.02$ s. The three-phase PWM excitation of period $T$ in the stator
under the no-load operation causes the $T$-periodic dynamics in
$\mathbf{u}$ which allows the imposition of the periodic constraint
\eqref{eq:dae_TP}. For more details regarding the mathematical model
we refer to \cite{Kulchytska-Ruchka_2018ac}. We would like to note
that equation \eqref{eq:dae} is a DAE of index-$1,$ which in case of
discretization with {Backward} Euler can be treated
essentially like an ODE~\cite{Schops_2018aa}.

\begin{figure}[t]
	\begin{subfigure}[b]{.59\textwidth}
		\centering
		\includegraphics[height=4.5cm]{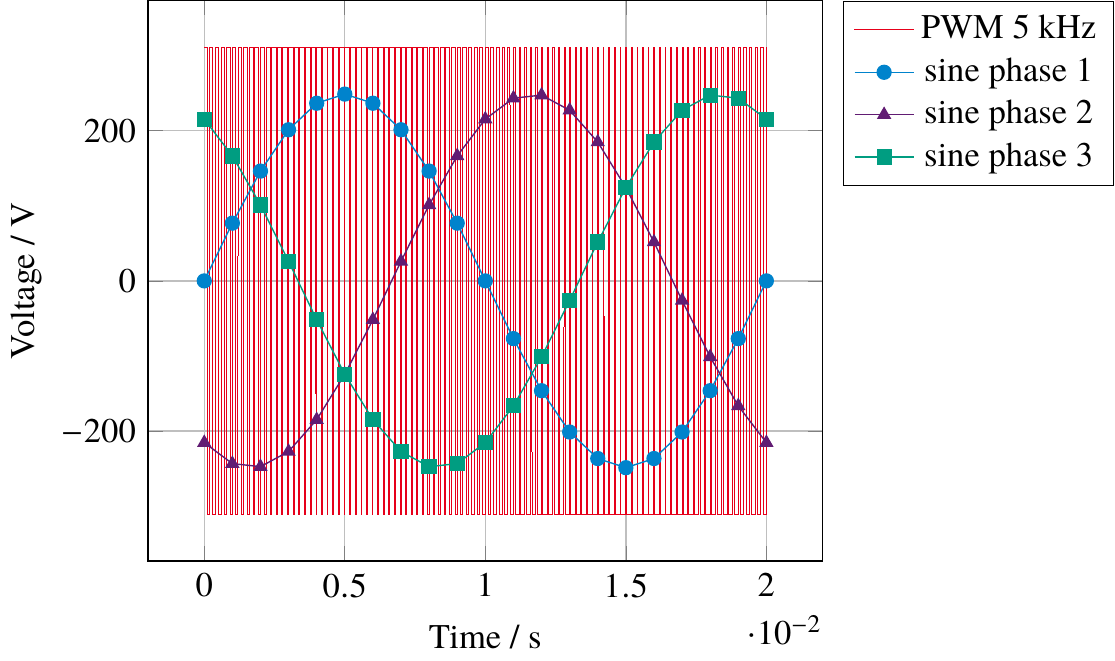}
	\end{subfigure}%
	\hspace*{0.05\textwidth}
	\begin{subfigure}[b]{.34\textwidth}
		\centering
		\includegraphics[height=4.5cm]{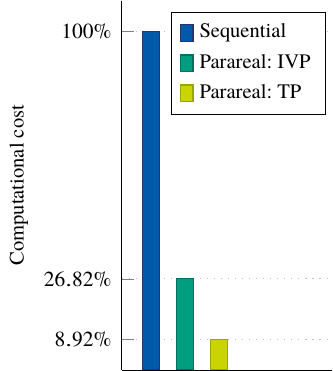}
	\end{subfigure}
  \caption{Left: PWM excitation \eqref{input_pwm} of $5$ kHz and
    three-phase sinusoidal voltage source of $50$ Hz, used as the
    coarse input in {our new periodic} Parareal algorithm
    \eqref{eq:PPIC_init}-\eqref{eq:PPIC_iter}.
    Right: comparison of the computational costs calculated in terms
    of the effective number of linear algebraic systems {solved} for
    different approaches to obtain the periodic steady-state solution
    of the induction machine model.}
  \label{fig:machine_comparison}
\end{figure}

{We now use our new periodic} Parareal {algorithm}
\eqref{eq:PPIC_init}-\eqref{eq:PPIC_iter} to find the solution of
\eqref{eq:dae}-\eqref{eq:dae_TP}. The fine propagator $\mathcal{F}$ is
then applied to \eqref{eq:dae} with the original three-phase PWM
excitation of $20$ kHz, discretized with the time step $\delta
T=10^{-6}$ s. The coarse solver $\bar{\mathcal{G}}$ uses a three-phase
sinusoidal voltage source with frequency $50$ Hz of the form
\eqref{input_sine}, discretized with the time step $\Delta T=10^{-3}$
s. Phase $1$ of the PWM signal switching at $5$ kHz as well as the
applied periodic coarse excitation on $[0,0.02]$ s are {shown} on
the left in Figure~\ref{fig:machine_comparison}.

Both coarse and fine propagators solve {the} IVPs for
\eqref{eq:dae} {using Backward} Euler, implemented within the GetDP
library \cite{Geuzaine_2007aa}.  We used $N=20$ {time} subintervals
for the simulation of the induction machine with {our new periodic
Parareal algorithm, which} converged after $14$ iterations. Within
these calculations $194\;038$ solutions of linearized systems of
equations were performed effectively, i.e, when considering 
the fine solution cost only on one subinterval 
(due to parallelization) together with the sequential coarse solves.

On the other hand, a classical way to obtain the periodic steady-state
solution is to apply a time integrator sequentially, starting from
{a} zero initial value at $t=0$. {This} computation reached the
steady state after $9$ periods, thereby requiring $2\;176\;179$ linear
system solves. Alternatively, one could apply the Parareal algorithm
with reduced coarse dynamics, introduced in
\cite{Kulchytska-Ruchka_2018ac}, to the IVP for \eqref{eq:dae} on
$[0,9T]$. In this case the simulation needed effectively
$583\;707$ sequential linear solutions due to parallelization. 
However, in practice one would not know the number of periods beforehand 
and one could not optimally distribute the time intervals. We visualize 
this data on the right {in} Figure~\ref{fig:machine_comparison}.
The{se results show that our new periodic Parareal algorithm}
\eqref{eq:PPIC_init}-\eqref{eq:PPIC_iter} with reduced coarse dynamics
directly delivers the periodic steady-state solution about $11$ times
faster than the standard time integration, and $3$ times {faster
than} the application of Parareal with the reduced 
dynamics to an IVP on $[0,9T]$.

\section{Conclusions}\label{section:conclusions}

{We introduced a new periodic Parareal algorithm with reduced
  dynamics}, which is able to efficiently handle quickly-switching
discontinuous excitations {in time-periodic problems. We
  investigated} its convergence {properties theoretically, and
  illustrated them} via application to a linear RL-circuit
example. {We then tested the performance of our new periodic
  Parareal algorithm in the} simulation of a two-dimensional model of
a four-pole squirrel-cage induction machine{, and a} significant
acceleration of convergence to the steady state {was} observed.  In
particular, with our new {periodic} Parareal algorithm {with
  reduced dynamics} it is possible to obtain the periodic solution
$11$ times faster than when performing the classical time stepping.
	
\bibliographystyle{spmpsci}
\bibliography{references}
	
\end{document}